\documentclass[11pt]{amsart}
\usepackage[utf8]{inputenc}

\usepackage{amsmath,amssymb}

\usepackage{color}

\def\cO{{\mathcal O}}

\newcommand{\ac}{\operatorname{ac}}

\newcommand{\GL}{\operatorname{GL}}

\newcommand{\Jac}{\operatorname{Jac}}
\newcommand{\res}{\operatorname{res}}

\let\sub=\subseteq

\def\eqc{{\mathrm{eqc}}}  
\def\heqc{h$^\eqc$}

\newcommand{\cL}{\mathcal{L}}

\newcommand{\ZZ}{\mathbb{Z}}
\newcommand{\FF}{\mathbb{F}}

\newcommand{\QQ}{\mathbb{Q}}
\newcommand{\RR}{\mathbb{R}}
\newcommand{\NN}{\mathbb{N}}

\newcommand{\Rbasic}{R^{\mathrm{basic}}}

\newcommand{\eqR}[1]{\overset{(R#1)}{=}}

\newtheorem{thm}{Theorem}[section]
\newtheorem{lem}[thm]{Lemma}
\newtheorem{lem-defn}[thm]{Lemma-Definition}
\newtheorem{cor}[thm]{Corollary}
\newtheorem{prop}[thm]{Proposition}

\theoremstyle{definition}
\newtheorem{defn}[thm]{Definition}
\newtheorem{notn}[thm]{Notation}

\newtheorem{rem}[thm]{Remark}
\newtheorem{conv-rem}[thm]{Convention-Remark}

\begin{document}

\author{Raf Cluckers}
\address{CNRS, Univ.~Lille,  UMR 8524 - Laboratoire Paul Painlev\'e, F-59000 Lille, France, and,
KU Leuven, Department of Mathematics, B-3001 Leu\-ven, Bel\-gium}
\email{Raf.Cluckers@univ-lille.fr}
\urladdr{http://rcluckers.perso.math.cnrs.fr/}

\author{Immanuel Halupczok}
\address{Lehrstuhl f\"ur Algebra und Zahlentheorie, Mathematisches Institut, Universit\"atsstr. 1, 40225 D\"usseldorf, Germany}
\email{math@karimmi.de}
\urladdr{http://www.immi.karimmi.de/en/}

\thanks{The authors would like to thank J.~Denef, F.~Loeser and A.~Macintyre for interesting discussions on the topics of the paper.
The author R.\,C. was partially supported by the European Research Council under the European Community's Seventh Framework Programme (FP7/2007-2013) with ERC Grant Agreement nr.\ 615722 MOTMELSUM, by KU Leuven IF C14/17/083, and thanks the Labex CEMPI  (ANR-11-LABX-0007-01). The author I.\,H. was partially supported by the \emph{SFB~878: Groups, Geometry and Actions}, by the research training group
\emph{GRK 2240: Algebro-Geometric Methods in Algebra, Arithmetic and Topology}, and by the individual research grant \emph{Archimedische und nicht-archimedische Stratifizierungen höherer Ordnung}, all three funded by the DFG.
Part of the work has been done while I.\,H. was affiliated to the University of Leeds.}

\dedicatory{Dedicated to Angus Macintyre, source of inspiration} 

\title[$p$-adic Kontsevich--Zagier integral operation rules]
{A $p$-adic variant of Kontsevich--Zagier integral operation rules and of Hrushovski--Kazhdan style motivic integration}



\begin{abstract}
We prove that if two semi-algebraic subsets of $\QQ_p^n$ have the same $p$-adic measure, then this equality can already be deduced using only some basic integral transformation rules. On the one hand, this can be considered as
a positive answer to a $p$-adic analogue of a question asked by Kontsevich--Zagier in the reals (though the question in the reals is much harder). On the other hand, our result can also be considered as stating that over $\QQ_p$, universal motivic integration (in the sense of Hrushovski--Kazhdan) is just $p$-adic integration.
\end{abstract}

\maketitle

\section{Introduction}

A period is a real number that can be obtained by integrating a rational function $f$ over a semi-algebraic domain $X \subset \RR^n$, where both $f$ and $X$ are defined with coefficients in $\QQ$.
Kontsevich--Zagier \cite{KZ.per} put forward the question whether manipulating real integrals by basic rules like Stokes, semi-algebraic change of variables, and linearity can explain all equalities between periods, possibly even in an algorithmic, decidable way.
On a different matter, Hrushovski--Kazhdan's \cite{HK.motInt}  version  of the motivic integral
(in algebraically closed valued fields of equi-characteristic $0$) is the universal map satisfying similar kinds of basic rules of integration; in other words, an equality between motivic integrals holds if and only if it follows from basic manipulations.

The main result of this paper is a positive answer to a $p$-adic analogue of the Kontsevich--Zagier question,  where one integrates
over semi-algebraic domains in $\QQ_p^n$ with respect to the $p$-adic measure: Any equality between such $p$-adic integrals can be deduced using a few specific basic manipulations. From the Hrushovski--Kazhdan point of view, this means that
the universal (motivic) integration theory over $\QQ_p$ is simply $p$-adic integration.
Note however that whereas the original question by Kontsevich--Zagier has deep links to transcendental number theory,
those links are completey lost in our $p$-adic case; notably, the $p$-adic measure of any $p$-adic semi-algebraic set is always a rational number.

Let us explain the $p$-adic variant in more detail. Suppose that $X_1 \sub \QQ_p^{n_1}, X_2 \sub \QQ_p^{n_2}$ are definable sets in the language of valued fields (also called semi-algebraic sets in the context of valued fields)
such that $\mu(X_i)$ is finite for $i = 1,2$, where $\mu$ denotes the Haar measure on $\QQ_p^{n_i}$, normalized so that $\ZZ_p^{n_i}$ has measure $1$.
Does $\mu(X_1) = \mu(X_2)$ imply that $X_1$ can be transformed into $X_2$ by some basic rules?
To make the question precise, we define a ring $R_{\{0\}}$ generated by sets $X_i$ as above, and quotient by relations corresponding to natural integral transformations (see Definition \ref{defn.pairs}).
A positive answer to the question then corresponds to: The map from this ring $R_{\{0\}}$ to $\RR$ sending the class of $X\subset \QQ_p^n$ to its measure $\mu(X)$ is injective; this is the statement of Corollary \ref{cor:bijection}.

Note that the ring $R_{\{0\}}$ is an analogue of the ring of values of universal motivic integration considered by Hrushovski--Kazhdan in \cite{HK.motInt} (though in our setting, some of the technicalities from \cite{HK.motInt} can be avoided). From that point of view, it is desirable to obtain a complete description of $R_{\{0\}}$, and indeed,
using some known results one rather easily obtains that the range of the above map is exactly $\QQ$, so in particular $R_{\{0\}}$ is isomorphic to $\QQ$.

The original question of Kontsevich--Zagier was not just about measures of sets $X$, but about integrals of functions on $X$. However, it is not very difficult to get from one version of the question to the other, using that integrals can be expressed in terms of measures of sets; see the explanations below Corollary~\ref{cor:bijection} for some details.

Our main result -- Theorem~\ref{thm.main} -- is a family version of Corollary \ref{cor:bijection}, namely: Given two definable families $(X_{i,s})_{s \in S}$ of sets, if $\mu(X_{1,s}) = \mu(X_{2,s})$ for each $s \in S$, then all transformations needed to turn $X_{1,s}$ into $X_{2,s}$ can be carried out uniformly in $s$ (i.e., with definable families of transformations). This is made precise by introducing a ring $R_S$ which is a family variant of the above ring $R_{\{0\}}$.
The proof builds on a similar kind of result in the value group, obtained in \cite{iC.grZ}, see Theorem \ref{thm.pres-fam} below.

One should also compare the basic transformation rules  (R1 -- R4)  and Corollary \ref{cor:bijection} with the classification of definable sets without measure from  \cite{Clu.DefQp}, where it is shown that two infinite definable subsets of $\QQ^n$ are definably isomorphic (namely, in definable bijection) if and only if they have the same dimension. The question to classify integrals rather than definable sets without measure was raized by Angus Macintyre in 2001, after \cite{Clu.DefQp}. We thank him for raizing this question.

\medskip 
We of course rely on classical model theory of $\cL$-definable sets enabled by the quantifier elimination result by Macintyre \cite{Mac.QE}, cell decomposition by Denef \cite{Den.cell}, and dimension theory by van den Dries \cite{Dri.dimDef}, to reduce to Presburger sets \cite{Pre.oagZ}. To go beyond the situation treated in this paper, we mention the upcoming analogue of the rings $R_S$ in elementary extensions of $\QQ_p$ in work in progress in the PhD thesis of Florian Severin. (Note that there is no Haar measure on non-standard elementary extensions.) Other generalizations, e.g.\ to the framework of motivic integrals from \cite{CL.mot} and finding adequate integral operation rules,  are left for the future. A more direct generalization to any finite field extension of $\QQ_p$ and to other languages than $\cL$ is formulated in the final Remark \ref{rem:genK-L}.

\section{Precise statement of the main results}
\label{sec.results}

Let $\cL = \{+,\cdot, \cO\}$ be the language of valued fields,
namely, with the ring operations and a predicate for the valuation ring.
By an $\cL$-definable set we mean a subset $X\subset \QQ_p^n$ for some $n$ which is given by a parameter free $\cL$-formula $\varphi$.  A function between $\cL$-definable sets is called $\cL$-definable if its graph is an $\cL$-definable set. (Remark~\ref{rem:genK-L} specifies some variants of the language for which our results also hold.)

When $X\subset \QQ_p^n$ is a Borel-measurable set, then we write $\mu(X) \in \RR_{\ge 0} \cup \{\infty\}$ for
the measure of $X$ with respect to the Haar measure 
on (the additive group) $\QQ_p^n$, normalized in such a way that $\mu(\ZZ_p^n) = 1$.
(This naturally also makes sense for $n = 0$, where $\QQ_p^0 = \ZZ_p^0$ is a one-point set, namely containig the empty tuple.)
It follows e.g.\ from Macintyre's quantifier elimination result \cite{Mac.QE} that any $\cL$-definable set is measurable.

\begin{defn}\label{defn.pairs}
Fix an $\cL$-definable set $S$.
We let $R_S$ be the abelian group generated by $\cL$-definable sets $X \sub S\times \QQ_p^n$ (for all $n \ge 0$) such that  $X_s:=\{x\in \QQ_p^n\mid (s,x)\in X\}$ has finite measure for each $s\in S$,  modulo the following relations, and where we write $[X]$ for the element of $R_S$ corresponding to the $\cL$-definable set $X$.

\begin{enumerate}
\renewcommand{\theenumi}{R\arabic{enumi}}
\item\label{r.add} (Additivity) If $X_1$ and $X_2$
are disjoint subsets of $S\times \QQ_p^n$, then
\[
[X_1 \cup X_2] = [X_1] + [X_2].
\]
\item\label{r.neg} (Negligable sets) If $X_s \sub \QQ_p^n$ has dimension less than $n$ for each $s$ in $S$, then
\[
[X] = 0.
\]
\item\label{r.cha} (Change of variables)
Suppose that $\phi\colon X \to Y$ is an $\cL$-definable bijection between $\cL$-definabe sets $X, Y \subset S \times \QQ_p^n$ inducing a bijection $\phi_s\colon X_s \to Y_s$ for each $s \in S$. Suppose moreover that
the sets $X_s$ and $Y_s$ are open in $\QQ_p^n$, that $\phi_s$ is $C^1$, and that the ($p$-adic) norm of the Jacobian determinant of $\phi_s$ equals $1$ everywhere, for each $s$ in $S$. Then
\[
[X] = [Y].
\]
\item\label{r.ball} (Product with unit ball) For any $\cL$-definabe $X \subset S \times \QQ_p^n$,
\[
[X] = [X \times \ZZ_p].
\]
\end{enumerate}
\end{defn}

The group $R_S$ can be endowed with the structure of commutative ring with unit, by defining $[X]\cdot [Y]$ as $[X \times_S Y]$, the class of the fiber product over $S$; see Lemma \ref{lem.product}.

For each fixed $s \in S$, the map sending
an $\cL$-definable set $X \sub S \times \QQ_p^n$ to the $p$-adic measure $\mu(X_s)$ factors over the relations
(R1)--(R4) (note that in (R3), $\mu(X_s) = \mu(Y_s)$ follows from the $p$-adic version of the change of variables formula for integrals). Therefore, this induces a map from $R_S$ to $\RR$, which clearly is a ring homomorphism:

\begin{defn}
Given a definable set $S$ and an element $s \in S$, we denote by $\mu_s\colon R_S \to \RR$ the (unique) ring homomorphism sending the class $[X]$ of an $\cL$-definable set $X \subset S \times \QQ_p^n$ to the $p$-adic measure $\mu(X_s)$ of its fiber at $s$.
\end{defn}

Now we can formulate our main result.

\begin{thm}\label{thm.main}
Fix an $\cL$-definable subset $S$ of some cartesian power of $\QQ_p$. Then the ring homomorphism
which sends an element $\Xi \in R_S$ to the function $S \to \RR, s \mapsto \mu_s(\Xi)$ is injective.
\end{thm}

It is well known (and can be deduced from from cell decomposition \cite{Den.cell}) that the $p$-adic measure of any definable set is a rational number. In particular,
the ring homomorphism from the above theorem takes values in the functions from $S$ to $\QQ$.
Describing its image precisely would be possible but rather technical. However, in the special case that $S$ is a singleton (e.g.~$S=\{0\}$), we can easily be more precise:

\begin{cor}\label{cor:bijection}
The map sending $[X]$ in $R_{\{0\}}$ to $\mu(X)$ is an isomorphism $R_{\{0\}} \to \QQ$ of rings.
\end{cor}

(The deduction of this corollary from the theorem is given after Remark~\ref{rem.balls}.)

As mentioned in the introduction, one can formulate a variant of the above results, where measures of sets are replaced by integrals of functions. One way to do this consists in considering a variant $R'_S$ of the above ring $R_S$ which is generated by pairs $(X,f)$, for $X \subset S \times \QQ_p^n$ and $f\colon X \to \QQ_p$ $\cL$-definable for which $\int_{X_s} |f(s, x)| d\mu(x)$ is finite for each $s$, where $|y| = p^{-v(y)}$ is the $p$-adic norm. The relations of $R'_S$ are natural analogues of the ones of $R_S$, where in (R3), we replace the assumption that the Jacobian determinant of $\phi_s$ has norm $1$ by the relation between the function on $X$ and the function on $Y$ coming from the change of variables formula.

To deduce the analogue of Theorem~\ref{thm.main} for $R'_S$, one uses that the natural map
from $R_S$ to $R'_S$ sending $[X]$ to $[X,1]$ is an isomorphism: Indeed, using the $R'_S$-version of (R3), one
obtains that an inverse can be defined by sending
$[X,f] \in R'_S$ to
$Y := \{(x,y) \in X \times \QQ_p \mid |y| \le |f(x)| \}$. (Note that $Y$ is defined in such a way that $\int_{X_s} f(s,\cdot)  = \mu(Y_s)$.)

\medskip

The main proof of Theorem~\ref{thm.main} is given in Section~\ref{sec.proof}; the strategy is as follows. Firstly, one uses Cell Decomposition to reduce to the case of certain definable sets $P(\Lambda)$ that are entirely described in terms of a Presburger definable set $\Lambda \sub \ZZ^n$. (Some variants of those $P(\Lambda)$ are also used.)
Then one further reduces to particularly simple such sets (called ``basic sets'' in Definition~\ref{defn.basic}). Those two reduction steps are carried out in the proof of Proposition~\ref{prop.make.basic.fam}.

The case of basic sets is treated in Proposition~\ref{prop.basic}: One reduces the problem of understanding equalities in $R_S$ to questions involving only Presburger definable sets, and those questions have already been answered in \cite{iC.grZ}.

The above proof requires a good understanding of the classes in $R_S$ of the sets $P(\Lambda)$ (and of their variants).
This understanding is developed in Section~\ref{sec.pres}.

\section{Sets defined in terms of Presburger data}
\label{sec.pres}

From now on, $S$ will almost always be a fixed $\cL$-definable subset of a cartesian power of $\QQ_p$, and $s$ will always be a variable running over $S$. We also use the following conventions:

\begin{conv-rem}
\begin{itemize}
\item
We will often identify an $\cL$-definable set $X \sub S\times \QQ_p^n$ with the family $(X_s)_{s \in S}$, and we call a family of subsets of $\QQ_p^n$ arising in this way an $\cL$-definable family (parametrized by $S$). Similarly, a family of maps $(f_s)_{s \in S}$ with $f_s\colon X_s \to Y_s$ is called an $\cL$-definable family if it arises from an $\cL$-definable map $f\colon X \to Y$, with $X \sub S \times \QQ_p^n$, $Y \sub S \times \QQ_p^m$.
 \item We consider the value group $\ZZ$ (together with $\infty$) as an imaginary sort, i.e., we call a subset of $\QQ_p^n \times (\ZZ \cup \{\infty\})^m$ an imaginary $\cL$-definable set if its preimage in
 $\QQ_p^{n+m}$ (under the map sending the last $m$ coordinates to their valuation) is $\cL$-definable. We may sometimes drop the word imaginary if it is clear from the context.
 (The notions of families are also applied in this generalized setting.)
 \item By Quantifier Elimination \cite{Mac.QE}, $X \subset \ZZ^n$ is an imaginary $\cL$-definable set if and only if it is a Presburger set, i.e., definable in the language $(+, <)$ of ordered abelian groups; we will sometimes also use this terminology.
\end{itemize}
\end{conv-rem}

\begin{notn}\label{notn.vtup}
For each integer $\ell> 0$, we write
$$
\ac_\ell:\QQ_p\to\ZZ_p/p^\ell\ZZ_p
$$
for the map sending $0$ to $0$ and nonzero $x$ to $xp^{-v( x)}\bmod p^\ell\ZZ_p$. We also write $\ac$ for $\ac_1$.
Given $x = (x_1, \dots, x_n) \in \QQ_p^n$, we write
\[
v(x) := (v(x_1), \dots, v(x_n)) \in \ZZ^n.
\]
\end{notn}

The relation (R3) in Definition~\ref{defn.pairs} is expresses that the class in $R_S$ of a definable set is preserved by a measure-preserving bijection. To formally work with (R3), we introduce the terminology
``(R3)-measure-preserving'' for maps $\phi$ satisfying the conditions of (R3); more precisely:

\begin{defn}
An $\cL$-definable map $\phi\colon X \to Y$ between $\cL$-definable sets
$X, Y \subset S \times \QQ_p^n$ is called (R3)-measure-preserving if for each $s \in S$, it induces
a bijection $\phi_s\colon X_s \to Y_s$ which is $C^1$ and such that the ($p$-adic) norm of the Jacobian determinant of $\phi_s$ is equal to $1$ everywhere. (For $C^1$ to make sense, we assume that $X_s$ is open.)
\end{defn}

Let us first state and prove that the group $R_S$ from Definition~\ref{defn.pairs} carries a natural ring structure.

\begin{lem}\label{lem.product}
Let $S$ be an $\cL$-definable set. The group $R_S$ (from Definition \ref{defn.pairs}) becomes a commutative ring with unit element $[S \times \ZZ_p]$ and with multiplication induced by
\[
[X]\cdot [Y] := [X \times_S Y],
\]
where $X \times_S Y$ is the fiber product over $S$, namely $\{(s,x,y)\mid s\in S,\ x\in X_s,\ y\in Y_s\}$.
\end{lem}

\begin{proof}
Let $\tilde R_S$ be the quotient of the free abelian group generated by same definable sets as in Definition~\ref{defn.pairs}, modulo the relation $[X] = [Y]$, whenever $X, Y \in S \times \QQ_p^n$ are definable sets such that there exists a definable bijection $\phi\colon X \to Y$ commuting with the projection to $S$. It is clear that
this group becomes a commutative ring using the above multiplication, and it is then equally easy to verify that that relations from Definition~\ref{defn.pairs} form an ideal in $\tilde R_S$, which finishes the proof.
\end{proof}

\begin{rem}\label{rem:inverse-inj}
Since $R_S$ is a ring with unit, we have a unique ring homomorphism $\ZZ \to R_S$. This ring homomorphism
is injective, since for any (fixed) $s$ in $S$, we obtain a one-sided inverse given by $[X] \mapsto \mu(X_s)$.
We will thus often identify $\ZZ$ with its image in $R_S$.
\end{rem}

\begin{rem}\label{rem:diff}
Any finite sum $\sum_i [X_i] \in R_S$ is equal to an element of the form $[Y]$. Indeed, using (\ref{r.ball}), we can first assume that $X_i \sub S \times \QQ_p^n \times \ZZ_p$ for all $i$ (where $n$ does not depend on $i$). Then we can translate the last coordinate of each $X_i$ to make them all  (using (R3)), and then we apply (\ref{r.add}).
In particular, every element of $R_S$ can be written in the form $[X] - [Y]$. Also note that one does not need to be careful concerning the order of coordinates (e.g.\ in (\ref{r.ball})), since (\ref{r.cha}) includes coordinate permutations.
\end{rem}

\begin{rem}\label{rem.S-union}
For disjoint $\cL$-definable sets $S, S'$, we have a natural
ring isomorphism $R_{S \cup S'} \cong R_S \times R_{S'}$.
\end{rem}

\begin{defn}\label{notn.presish}
Given an imaginary $\cL$-definable set $\Lambda \sub S \times \ZZ^n$ (for some $n \ge 0$), we define a set $P(\Lambda) \subset S \times \QQ_p^n$ via its fibers over $S$:
\[
P(\Lambda)_s =
\{x = (x_1, \dots, x_n) \in \QQ_p^n \mid \ac(x_1) = \dots =
\ac(x_n) = 1, v(x) \in \Lambda_s\}
\]
for $s \in S$.

If, in addition, we are given an imaginary $\cL$-definable function $\nu$ to $\ZZ$ whose domain contains $\Lambda$, we similarly define
$P(\Lambda, \nu) \subset S \times \QQ_p^n$ via
\[
P(\Lambda, \nu)_s =
\{(x, y) \in P(\Lambda)_s \times \QQ_p \mid \ac(y)=1,\
v(y) =  - n -1 - \nu_s(v(x)) -\sum_{i=1}^n v(x_i)
\}.
\]
\end{defn}

Clearly, the sets $P(\Lambda)$ and $P(\Lambda, \nu)$ are definable.
The motivation behind this definition of $P(\Lambda, \nu)$ is that it is a simple way of defining a set with a prescribed $p$-adic measure, namely (as an easy computation shows)
\begin{equation}\label{eq.bal.reason}
\mu(P(\Lambda, \nu))_s = \sum_{\lambda \in \Lambda_s} p^{\nu_s(\lambda)}.
\end{equation}

We will see, in several of the following lemmas, that the classes in $R_S$ of sets of the form $P(\Lambda, \nu)$
satisfy relations one would expect from (\ref{eq.bal.reason}) (provided that $\mu(P(\Lambda, \nu)_s)$ is finite for all $s$).

\begin{lem}\label{lem.pn-1}
Let $\Delta_n$ be the image of the diagonal embedding of the non-negative integers $\NN$ into $\ZZ^n$, i.e.,
\[
\Delta_n = \{(\underbrace{\lambda, \dots, \lambda}_{n}) \in \ZZ^n  \mid \lambda \in\NN\}.
\]
Then we have $(p^n-1)\cdot [P(S \times \Delta_n)] = 1$ in the ring $R_S$.
\end{lem}

\begin{proof}
For any $a \in \FF_p^n \setminus \{0\}$, there exists a matrix $M \in \GL_n(\ZZ_p)$ such that $\res(M)$ sends the vector $(1, \dots, 1) \in \FF_p^n$ to $a$.
Such an $M$ is (R3)-measure-preserving, and it sends $P(\Delta_n)$ to
\[
X_{a} := \{x\in \ZZ_p^n\mid  \forall i,j\colon v(x_i)=v(x_j) \text{ and } \forall i\colon \ac(x_i)=a_i  \}.
\]
By (R3), we obtain $[S \times X_{a}] = [S \times P(\Delta_n)]$, and
since the ($p^n - 1$ many) sets $X_a$ (for $a$ as above) form a partition of $\ZZ_p^n \setminus \{0\}$, we obtain (in $R_S$):
\[(p^n-1)[P(S \times \Delta_n)] \eqR1 [S \times (\ZZ_p^n \setminus \{0\})] \eqR2 [S \times \ZZ_p^n] \eqR4 1.\qedhere\]
\end{proof}

\begin{lem}\label{lem.ac-change}
Given an imaginary $\cL$-definable subset $\Lambda \sub S \times \ZZ$ and an integer $\ell >0$,
define $P_\ell(\Lambda)$ via
$P_\ell(\Lambda)_s = \{x \in P(\Lambda)_s \mid \ac_\ell(x) = 1\}$. Then we have
$p^{\ell - 1}[P_\ell(\Lambda)] = [P(\Lambda)]$ in $R_S$. As similar statement holds
for $P_\ell(\Lambda, \nu)_s := \{x \in P(\Lambda,\nu)_s \mid \ac_\ell(x) = 1\}$, when additionally
an $\cL$-definable $\nu\colon \Lambda \to \ZZ$ is given.
\end{lem}

\begin{proof}
Choose representatives $r_1, \dots, r_{p^{\ell-1}}$ of the different cosets of $(1 + p \ZZ_p) / (1 + p^\ell \ZZ_p)$.
Then multiplication by $r_i$ is (R3)-measure-preserving and for each $s$, the sets $r_i P_\ell(\Lambda)_s$
form a partition of $P(\Lambda)_s$. Thus the claim follows by (R3) and (R1).

The same proof also gives the second part, if in ``$r_i P_\ell(\Lambda,\nu)_s$'', one lets $r_i$ act on the first coordinate only.
\end{proof}

\begin{lem}\label{lem.div}
There exists a (unique) injective ring homomorphism $\QQ \to R_S$. In particular, the additive group of $R_S$ is divisible and torsion free.
\end{lem}

\begin{proof}
By Remark \ref{rem:inverse-inj}, we have $\ZZ \subset R_S$, so it suffices to prove that $R_S$ contains a multiplicative inverse of $\ell$ for every prime $\ell$. (Indeed, torsion freeness then follows by multiplying both sides of an equation of the form ``$n\cdot \Xi = 0$'' (where $n \ge 1$ and $\Xi \in R_S$) by $\frac1n$.)

Since $\ZZ_p$ is the disjoint union of $p$ translates of $p\ZZ_p$, we have $p\cdot [S \times p\ZZ_p] = 1$, so that $[S \times p\ZZ_p]$ is a multiplicative inverse of $p$.
If $\ell \ne p$, then Lemma~\ref{lem.pn-1} provides the desired multiplicative inverse (namely, a multiple of $[P(S \times \Delta_n)]$), provided that we can find an $n \ge 1$ such that $\ell$ divides $p^n-1$. Indeed,
the image of $p$ in the ring $\ZZ/\ell\ZZ$ is a unit, so for $n$ the order of that image
in the group $(\ZZ/\ell\ZZ)^\times$ of units, we obtain $p^n \equiv 1 \mod \ell$ and hence $\ell$ divides $p^n-1$.
\end{proof}

We will from now on identify $\QQ$ with its image in $R_S$.

\begin{rem}\label{rem.balls}
We now easily see that the class in $R_S$ of a ball in $\QQ_p$ is the expected one, namely:
The above proof in particular yields that $[S \times p\ZZ_p] = \frac1p$. In a similar way,
we deduce $[S \times p^c\ZZ_p] = p^{-c}$ for any integer $c$, and
then also $[S \times (a+ p^c\ZZ_p)] = p^{-c}$ for any $a \in \QQ_p$ (by applying (R3) to the translation by $a$).
\end{rem}

We now have all the ingredients to deduce Corollary~\ref{cor:bijection} from Theorem~\ref{thm.main}:

\begin{proof}[Proof of Corollary~\ref{cor:bijection} using Theorem~\ref{thm.main}]
By Theorem~\ref{thm.main}, the map from $R_{\{0\}}$ to $\QQ$ induced by the Haar measure is injective. Its restriction to $\QQ \sub R_{\{0\}}$ is clearly the identity, so surjetivity follows.
\end{proof}

\begin{lem}\label{lem.pres-prod}
For $i = 1,2$ and $n_i \in \NN$, let
$\Lambda_i \sub S \times \ZZ^{n_i}$
be imaginary $\cL$-definable sets, let $\nu_{i}\colon \Lambda_{i} \to \ZZ$ be imaginary
$\cL$-definable functions,
and define $\nu\colon \Lambda_1 \times_S \Lambda_2 \to \ZZ$ by $(s, \lambda_1,\lambda_2) \mapsto \nu_{1}(s,\lambda_1) + \nu_{2}(s,\lambda_2)$.
Suppose that $\Lambda_{i,s}$ is finite for $i=1,2$ and each $s$ in $S$.
Then we have
\[
[P(\Lambda_{1} \times_S \Lambda_{2},\nu)]
=
[P(\Lambda_{1}, \nu_{1})]\cdot
[P(\Lambda_{2}, \nu_{2})]
\]
in $R_S$.
\end{lem}

\begin{rem}\label{rem.pres-prod}
Some useful special cases are obtained for $n_2 = 0$: Let $\Lambda \sub S \times \ZZ^n$ and $\nu\colon \Lambda \to \ZZ$ be $\cL$-definable, with $\mu(P(\Lambda, \nu)_s)$ finite for every $s \in S$.
\begin{enumerate}
 \item If $\nu$ is of the form $\nu(s,\lambda) = \nu'(s)$ for some $\nu'\colon S \to \ZZ$, then
   $[P(\Lambda, \nu)] = [P(\Lambda \times_S S, \nu)] = [P(\Lambda, 0)] \cdot [P(S, \nu')]$.
 \item If $c$ is an integer, then
   $[P(\Lambda, \nu+c)] = [P(\Lambda \times_S S, \nu+c)] = [P(\Lambda, \nu)] \cdot [P(S, c)] = p^c\cdot [P(\Lambda, \nu)]$,
   where the last equality holds by Remark~\ref{rem.balls}.
\end{enumerate}
\end{rem}

\begin{proof}[Proof of Lemma~\ref{lem.pres-prod}]
To prove the lemma, we will specify an $\cL$-definable family of (R3)-measure-preserving
bijections
\[P(\Lambda_{1} \times_S \Lambda_{2},\nu)_s \times \ZZ_p \to P(\Lambda_{1}, \nu_{1})_s \times
P(\Lambda_{2}, \nu_{2})_s.\]
The lemma then follows from (R3), (R4), and the definition of the product.

Fix $s \in S$. We have two natural projections:
\[
\pi_{L}\colon P(\Lambda_{1} \times_S \Lambda_{2}, \nu)_s \to P(\Lambda_{1} \times_S \Lambda_{2})_s
\]
and
\[
\pi_R\colon P(\Lambda_{1}, \nu_{1})_s \times
P(\Lambda_{2}, \nu_{2})_s \to P(\Lambda_{1})_s \times
P(\Lambda_{2})_s =  P(\Lambda_{1} \times_S \Lambda_{2})_s
\]
For $x \in P(\Lambda_{1} \times_S \Lambda_{2})_s$,
the fibers over $x$ are of the form
$$
\pi_{R}^{-1}(x) = p^{\sigma}+p^{\sigma+1}\ZZ_p
$$
and
$$
\pi_{L}^{-1}(x) = (p^{\sigma_1}+p^{\sigma_1+1}\ZZ_p) \times
(p^{\sigma_2}+p^{\sigma_2+1}\ZZ_p)
$$
for some $\sigma, \sigma_1, \sigma_2$ depending definably on $x$ and $s$ and satisfying
$\sigma = \sigma_1 + \sigma_2+1$.
Choose an $\cL$-definable family of functions $h_s$ sending $x$ in $P(\Lambda_{1} \times_S \Lambda_{2})_s$ to an element of $\QQ_p$ of valuation $\sigma_2+1$ and with $\ac(h_s(x))=1$; such an $\cL$-definable family exists by existence of definable Skolem functions; see \cite{Dri.skolQp}.
Using this, we obtain an (R3)-measure-preserving bijection
from $\pi_{R}^{-1}(x)\times \ZZ_p$ to $\pi_{L}^{-1}(x)$ sending $(a,b)$ to $(a/h_s(x) , (1/p+ b) h_s(x) )$.
Gluing those bijections together for all $x$ yields the desired (R3)-measure-preserving bijection.
\end{proof}

The following proposition states that
a definable bijection $\Lambda_1 \to \Lambda_2$ at the value group level induces equality
between the classes $[P(\Lambda_i)]$ in $R_S$.

\begin{prop}\label{prop.Z-to-Qp}
For $i = 1,2$ and $n_i \in \NN$, let
$\Lambda_{i} \sub S \times \ZZ^{n_i}$
be imaginary $\cL$-definable sets and let $\nu_{i}\colon \Lambda_{i} \to \ZZ$ be
imaginary $\cL$-definable functions. Suppose that the sets $P(\Lambda_{i}, \nu_{i})_s$ have finite measure for each $i$ and each $s$.
Suppose moreover that there exists an imaginary
$\cL$-definable family of bijections
$\phi_s\colon \Lambda_{1,s} \to \Lambda_{2,s}$ which is compatible with the $\nu_{i}$, i.e., such that $\nu_{1}(s, \lambda) = \nu_{2}(s, \phi_s(\lambda))$ for every $(s,\lambda) \in \Lambda_{1}$. Then we have
$$
[P(\Lambda_{1}, \nu_{1})] = [P(\Lambda_{2}, \nu_{2})]
$$
in $R_S$.
\end{prop}

The idea of the proof is to reduce to some very simply cases, by decomposing $\Lambda_i$ into finitely many pieces and by writing $\phi_s$ as a composition of finitely many maps. In those simple cases, explicit definable measure-preserving bijections can be obtained using the existence of definable Skolem functions.
Here are the details:

\begin{proof}[Proof of Proposition~\ref{prop.Z-to-Qp}]
We may suppose $n_1 = n_2$, since if, say, $n_1 < n_2$,
we can replace $P(\Lambda_{1}, \nu_{1})$ by
\[
P(\{-1\}^{n_2-n_1} \times \Lambda_{1}, \nu_{1} \circ \pi_{>n_2-n_1})
=
P(\{-1\})^{n_2-n_1} \times P(\Lambda_{1}, \nu_{1}),
\]
where $\pi_{>n_2-n_1}\colon \QQ_p^{n_2} \to \QQ_p^{n_1}$ is the projection to the last $n_1$ coordinates.
Note that $P(\{-1\})$ is just a translate of $\ZZ_p$.

Now that $n_1 = n_2 =: n$,
one technique used in this proof consists in constructing
an $\cL$-definable family of $C^1$ bijections $\psi_s\colon P(\Lambda_{1})_s \to P(\Lambda_{2})_s$ such that, for every $s \in S$ and every $\lambda \in \Lambda_{1,s}$, $\psi_s$ restricts to a bijection $P(\{\lambda\}) \to P(\{\phi_s(\lambda)\})$ and such that $v(\Jac\psi_s (x))$ is constant (and finite) on $P(\{\lambda\})$. This then induces a bijection
\[
\tilde{\psi}_s\colon P(\Lambda_{1})_s \times \QQ_p \to P(\Lambda_{2})_s \times \QQ_p,
(x,y) \mapsto (\psi(x), y/ J_s(x)),
\]
where $J_s$ is a chosen $\cL$-definable family of $C^1$-functions with $\ac(J_s(x))=1$ and $v(J_s(x))=v(\Jac\psi_s(x))$ for all $x$ in $P(\Lambda_{1})_s$ (such $J$ exists by existence of Skolem functions; see \cite{Dri.skolQp}).
An easy computation (using the compatibility of $\phi_s$ with the $\nu_i$) shows that $\tilde{\psi}_s$ restricts to a bijection
$P(\Lambda_{1}, \nu_{1})_s \to P(\Lambda_{2}, \nu_{2})_s$ and that this restriction satisfies (R3). Thus, the proposition follows whenever we can find
$\psi_s$ as above.

Instead of applying this technique directly in general, we first reduce to special cases where the given family $(\phi_s)_{s \in S}$ is of a simple form.

Given a partition of $\Lambda_1$ into finitely many $\cL$-definable sets $\Lambda'_1$, it suffices to prove the lemma
for those $\Lambda'_1$ and $\Lambda'_2$ given by $\Lambda'_{2,s} = \phi_s(\Lambda'_1)$. (Indeed, such a partition induces corresponding partitions of $P(\Lambda_i, \nu_i)$; then use (R1).)
By piecewise linearity of Presburger functions, using such a finite partition, we may assume that $\phi_s$ is of the form
\[
\phi_s(\lambda) = M\lambda + \mu_s
\]
for some matrix $M = (m_{ij})_{ij}$ and some vector $\mu_s = (\mu_{s,i})_i$, both with coefficients in $\QQ$, and where $M$ does not depend on $s$. Since $\phi_s$ is a bijection, we may moreover assume that $M$ is invertible.
(Note that this might require a refinement of the finite partition.)

If we can write $\phi_s$ as a composition of several maps $\phi_{k,s}$ (each one forming an $\cL$-definable family),
it suffices to prove the proposition for each of the $\phi_{k,s}$. In this way, we may further assume that $\phi_s$ is of one of the following forms:
\begin{enumerate}
\item $\lambda \mapsto \lambda + \mu_s$ for some $\mu_s \in \ZZ^n$ (where $s \mapsto \mu_s$ is $\cL$-definable);
\item a permutation of coordinates;
\item $(\lambda_1, \dots, \lambda_n) \mapsto (\lambda_1+ \lambda_2,\,\, \lambda_2, \dots, \lambda_n)$.
\item $(\lambda_1, \dots, \lambda_n) \mapsto (r\lambda_1,\,\, \lambda_2, \dots, \lambda_n)$ for some $r \in \QQ$, $r \ne 0$;
\end{enumerate}
(In (2)--(4), everything is independent of $s$.)
We now prove the proposition in each of these cases, partly by specifying a family $\psi_s$ as required for the technique described at the beginning of the proof.

(1)
By \cite{Dri.skolQp}, there exists an $\cL$-definable function sending $s \in S$ to an element $(a_1, \dots, a_n) \in \QQ_p^n$ satisfying $v(a_i) = \mu_{s,i}$ and $\ac(a_i)=1$ for each $i$. Then we define $\psi_s(x_1, \dots, x_n) = (a_1x_1, \dots, a_nx_n)$. As required, this sends $P(\{\lambda\})$ to $P(\{\lambda + \mu_s\})$, and the valuation of its Jacobian is constant (namely equal to $\sum_i \mu_i$).

(2)
Clear.

(3)
Set $\psi_s(x_1, \dots, x_n) := (x_1\cdot x_2, \, x_2, \dots, x_n)$.

(4)
We may assume $r \in \ZZ$, since then we can obtain arbitrary $r$ by composing one such map with an inverse of such a map. Also, without loss we may suppose that $n = 1$.

By Hensel's Lemma one easily finds (see \cite[Corollary 1]{Clu.DefQp}) that there exists $\ell>0$ such that the map
$$
\psi\colon x \mapsto x^r
$$
defines a bijection
\[
(\ac_{\ell})^{-1}(1) \to (\ac_{\ell'})^{-1}(1)\cap P(r\ZZ),
\]
where
$$
\ell' = \ell + v(r).
$$
Using the notation from Lemma~\ref{lem.ac-change},
$\psi$ restricts to bijections $P_\ell(\{\lambda\}) \to P_{\ell'}(\{r\lambda\})$ for each $\lambda$,
and for each such restriction, the valuation of the Jacobian is constant, namely: for $x \in  P_\ell(\{\lambda\})$,
that valuation is equal to $v(rx^{r-1}) = v(r) + (r-1)\cdot \lambda$. By a variant of the technique from the beginning
of the proof, we obtain an $\cL$-definable bijection
\[
\tilde\psi\colon P_\ell(\Lambda_1, \nu_1) \to P_{\ell'}(\Lambda_2, \nu_2)
\]
whose Jacobian has valuation constant equal to $v(r)$. From this, one deduces that $[P_{\ell'}(\Lambda_2, \nu_2)] = p^{-v(r)}\cdot [P_\ell(\Lambda_1, \nu_1)]$. Now the proposition in Case (4) follows using Lemma~\ref{lem.ac-change}.
\end{proof}

We end this section with two results from \cite{iC.grZ} about Presbuger definable families that will be needed in the proof of Theorem~\ref{thm.main}.

To apply Proposition~\ref{prop.Z-to-Qp}, we need to find definable families $\phi_s\colon \Lambda_{1,s} \to \Lambda_{2,s}$
in the value group. Those will be obtained using the following variant of Theorem~\ref{thm.main} for Presburger definable sets:

\begin{thm}[{\cite[Theorem~5.2.2]{iC.grZ}}]\label{thm.pres-fam}
Let $\tilde S \subseteq \ZZ^k$ be a Presburger set.
Let $\Lambda_{\tilde s}$ and $\Lambda'_{\tilde s}$ be two Presburger families, where $\tilde s$ runs over $\tilde S$.
Suppose moreover that for each $\tilde s \in \tilde S$, $\Lambda_{\tilde s}$ and $\Lambda'_{\tilde s}$ are finite sets of the same cardinality.
Then there exists a Presburger family of bijections $\phi_{\tilde s}\colon \Lambda_{\tilde s} \to \Lambda'_{\tilde s}$, with $\tilde s$ running over $\tilde S$.
\end{thm}

(Recall that Presburger set is just an imaginary $\cL$-definable subset of $\ZZ^n$, and similarly for Presburger maps.)

In the proof of Theorem~\ref{thm.main}, we will also need to understand how, in the above setting,
the cardinality of $\Lambda_{\tilde s}$ can depend on $\tilde s$. The following proposition states that the dependence is piecewise polynomial:

\begin{prop}[{\cite[Proposition~5.2.1]{iC.grZ}}]\label{prop.pres-pol}
Let $\tilde S \subseteq \ZZ^k$ be a Presburger set, and
let $\Lambda_{\tilde s}$ be a Presburger family of finite sets, where $\tilde s$ runs over $\tilde S$.
Then there exists a partition of $\tilde S$ into finitely many Presburger sets $\tilde S_i$
and polynomials $g_{i} \in \QQ[x_1, \dots, x_k]$ such that $\# X_{\tilde s} = g_i(\tilde s)$ for each $\tilde s\in \tilde S_i$.
\end{prop}

\section{Proof of Theorem \ref{thm.main}}
\label{sec.proof}

We now have all the ingredients to prove our main result, Theorem~\ref{thm.main}, namely that an element $\Xi$ of $R_S$ is determined by the measures $\mu_s(\Xi)$, for $s \in S$.
For that proof, we first introduce a subring of $R_S$, generated by sets that are unions of finitely many boxes, all of which have the same measure, but where the number of boxes may depend on $s \in S$. We will then first prove  Theorem~\ref{thm.main} for that subring (Proposition~\ref{prop.basic}), and then show that $R_S$ is actually not much bigger than the subring (Proposition~\ref{prop.make.basic.fam}); this will then easily imply the theorem.

\begin{defn}\label{defn.basic}
We let $\Rbasic_S$ be the subgroup of $R_S$ generated by the classes $[P(\Lambda, \nu \circ \pi_S)]$  of
sets of the form $P(\Lambda, \nu \circ \pi_S)$ for $\Lambda \sub S\times \ZZ^n$ and $\nu \colon S \to \ZZ$ imaginary $\cL$-definable (for some $n \ge 0$), where $\pi_S: \Lambda \to S$ is the projection to the $S$-coordinates, and where moreover
$\Lambda_s$ is finite for every $s \in S$. We call such $P(\Lambda, \nu \circ \pi_S)$ a \emph{basic} set,
and by abuse of notation, we also denote it by $P(\Lambda, \nu)$ (thinking of $\nu$ as a function on $\Lambda$ only depending on the $S$-variable).
\end{defn}

\begin{rem}
By Lemma \ref{lem.pres-prod},  $\Rbasic_S$ is a subring of $R_S$.
\end{rem}

Now we are ready to do the first main step of the proof of Theorem~\ref{thm.main}.

\begin{prop}\label{prop.basic}
Suppose that $\Xi \in \Rbasic_S$ is an element
satisfying $\mu_s(\Xi) = 0$ for every $s \in S$.
Then $\Xi = 0$.
\end{prop}

\begin{proof}
Let $\Xi \in \Rbasic_S$ be given, satisfying $\mu_s(\Xi) = 0$ for all $s$ in $S$. (We need to prove $\Xi = 0$.) Write $\Xi$ as a finite sum of generators of $\Rbasic_S$, i.e.,
\[\Xi = \sum_j \delta_j \cdot  [P(\Lambda_j,\nu_j)]
 \]
with the $P(\Lambda_j,\nu_j)$ basic sets and
with $\delta_j$ either $1$ or $-1$.

Given a partition of $S$ into finitely many $\cL$-definable sets $S_i \sub S$, we get natural images $\Xi_i$ of $\Xi$ in $R_{S_i}$ (by Remark~\ref{rem.S-union}), and it suffices to prove that $\Xi_i = 0$ in $R_{S_i}$ for each $i$.
This allows us to apply cell decomposition to $S$, so that we may without loss assume that
$S$ equals $P(\tilde{S})$ for some Presburger set $\tilde{S} \sub \ZZ^m$ (for some $m$).
Using that the sets $\Lambda_{j,s}$ and the values $\nu_{j}(s)$ live in the value group,
we may moreover assume (by choosing the cell decomposition appropriately) that
they factor through the coordinate-wise valuation map $v:S\to \tilde{S}$, i.e.,  $\Lambda_{j,s_1} = \Lambda_{j,s_2}$ and $\nu_j(s_1) = \nu_j(s_2)$ whenever $v(s_1) = v(s_2)$.

We write $\Lambda_{j,\tilde s}$ and $\nu_j(\tilde s)$ for $\Lambda_{j,s}$ and $\nu_j(s)$, respectively,
when $s\in S$ and $\tilde s\in \tilde{S}$ satisfy $v(s)=\tilde s$. We denote the coordinates of $\tilde s$ by $\tilde s_i$ ($i = 1, \dots, m$).

By Proposition~\ref{prop.pres-pol} (and using Remark~\ref{rem.S-union} once more to partition $\tilde{S}$), we may assume that
$$
\#\Lambda_{j,\tilde{s}} = g_j(\tilde{s})
$$
for some polynomials $g_j \in \QQ[\tilde{s}_1, \dots, \tilde{s}_m]$
and that $\nu_j$ is affine linear with coefficients in $\QQ$, i.e., of the form
\begin{equation}\label{eq.c.b}
\nu_j(\tilde{s}) = c_j + \underbrace{\sum_i b_{j,i}\cdot \tilde{s}_i}_{=: b_j\cdot \tilde s}
\end{equation}
for some rational numbers $c_j$ and $b_{j,i}$. Using this notation,
we obtain, for $s \in S$ and $\tilde s = v(s)$:
\[
\mu_s(\Xi) = \sum_j \delta_j \mu(P[\Lambda_j, \nu_j])
= \sum_j \delta_j \cdot g_j(\tilde s)\cdot p^{\nu_j(\tilde s)},
\]
so the assumption that $\mu(\Xi_s) = 0$ for all $s$ becomes
\begin{equation}\label{eq.sum-poly}
\sum_j \delta_j \cdot g_j(\tilde{s})\cdot p^{\nu_j(\tilde{s})} = 0 \qquad
\text{for all } \tilde{s} \in \tilde{S}.
\end{equation}
By Rectilinearization  \cite[Theorem 2]{Clu.cell}, Proposition \ref{prop.Z-to-Qp}
and a further finite partition of $\tilde{S}$, we may then assume that
$\tilde{S} = \NN^m$. After this modification, we in particular get that the $g_j$ have integer coefficients and that the $c_j$ and $b_{i,j}$ are integers.

We may assume that the constant terms $c_j$ from (\ref{eq.c.b}) are non-negative; if not, we replace $\Xi$ by
\[
\Xi' := \sum_j \delta_j\cdot [P(\Lambda_j, \nu_j + c)]
\]
for a suitable $c$. (Note that by Remark~\ref{rem.pres-prod}, we have $\Xi' = p^c \cdot\Xi$ in $R_S$, so $\Xi'=0$ implies $\Xi = 0$, by torsion freeness of $R_S$.)

Now we can entirely get rid of the $c_j$: Intuitively, replacing $\nu_j$ by $\nu_j - c_j$ divides $\Xi$ by $p^{c_j}$; we make up for this by replacing $\Lambda_j$ by $p^{c_j}$ disjoint copies of itself. Formally, this is the following computation
(which uses Lemma \ref{lem.pres-prod} multiple times):
\begin{align*}
[P(\Lambda_j,\nu_j)]
&=
[P(\Lambda_j ,\nu_j - c_j)] \cdot [P(S, c_j)]
=
[P(\Lambda_j, \nu_j - c_j)] \cdot p^{c_j}\\
&=
[P(\Lambda_j, \nu_j - c_j)] \cdot [P(\Lambda'_j, 0)]
=
[P(\Lambda_j \times \Lambda'_j, \nu_j - c_j)],
\end{align*}
where $\Lambda'_j$ a Presburger set with $p^{c_j}$ elements.

Now that the $c_j$ are gone, we may group summands of (\ref{eq.sum-poly}) together that have the same sign $\delta_j$ and the same tuple $(b_{j,1}, \dots, b_{j,m})$. Here, ``grouping'' summands $j$ and $j'$ means first making $\Delta_j$ and $\Delta_j'$
disjoint (using definable bijections) and then taking their union. In this way, (\ref{eq.sum-poly}) becomes, after some relabeling:
 \begin{equation}\label{eq.sum-both}
\sum_{j} (g_{j^+}(\tilde{s}) -  g_{j^-}(\tilde{s}))\cdot p^{b_j \cdot \tilde{s}} = 0 \qquad
\text{for all } \tilde{s} \in \tilde{S} = \NN^m,
\end{equation}
where $b_j \ne b_{j'}$ for $j \ne j'$.
Using a suitable induction according to growth rates of the sum (as a function of the $\tilde s_i$), (\ref{eq.sum-both}) implies $g_{j^+} = g_{j^-}$ for all $j$. In particular, for every $\tilde{s} \in \tilde{S}$, the corresponding sets $\Lambda_{j^+,\tilde s}$ and $\Lambda_{j^-,\tilde s}$ have the same cardinality.
By Theorem \ref{thm.pres-fam}, this implies that there exists an $\cL$-definable family of bijections
$\phi_s\colon \Lambda_{j^+,\tilde s} \to \Lambda_{j^-,\tilde s}$. Then Proposition \ref{prop.Z-to-Qp} yields
$[P(\Lambda_{j^+}, \nu_j)] = [P(\Lambda_{j^-}, \nu_j)]$, which, by summing over $j$, implies $\Xi = 0$.
\end{proof}

\begin{prop}\label{prop.make.basic.fam}
For every $\Xi \in R_S$, there exists an integer $\ell>0$ such that $\ell\cdot \Xi \in \Rbasic_S$.
\end{prop}

\begin{proof}
Since every element of $R_S$ can be written as a difference of two generators (by Remark~\ref{rem:diff}), it suffices to deal with the case $\Xi = [X]$, for some $\cL$-definable $X \subset S \times \QQ_p^n$ (satisfying that $X_s$ has finite $p$-adic measure for each $s \in S$).

By cell decomposition (and since $\Rbasic_S$ is closed under finite sums), we may assume that $X$ is a cell over $S$ in the sense of e.g.\ \cite[Definition~3.4]{CCL.lipschitz} (which ultimately originates from \cite{Den.rat}, proof of Theorem 7.4).
We may also assume that $\dim X = n$, since $[X] = 0$ for lower-dimensional $X$.
Moreover, we may get rid of the centers of the cell by (R3)-measure-preserving translations. In this way, we reduce to the case where $X$ is of the form
\[
X = \{(s,x) \in S \times \QQ_p^n \mid v(x) \in \Lambda_s, \ac_{\ell}(x_1) = \xi_1, \dots, \ac_{\ell}(x_n) = \xi_n
\},
\]
for an imaginary $\cL$-definable set $\Lambda \sub S \times \ZZ^n$, an integer $\ell>0$ and some $\xi_i \in \ZZ_p/p^{\ell}\ZZ_p$ (where $x = (x_1, \dots, x_n)$ and $v(x) = (v(x_1), \dots, v(x_n))$).
Moreover, we can replace all $\xi_i$ by $1$ (by multiplying each coordinate by a suitable element of $\ZZ_p^\times$),
so that $X = P_\ell(\Lambda)$, in the notation of Lemma~\ref{lem.ac-change}. Finally,
that lemma allows us to reduce to the case where $\ell=1$ (so that $X = P(\Lambda)$).

Next, we replace $X$ by $P(\Lambda, \nu)$, where $\nu\colon \Lambda \to \ZZ$ is chosen in such a way that $P(\Lambda, \nu) = P(\Lambda) \times (p^{-1} + \ZZ_p)$ (which implies $[P(\Lambda)] = [P(\Lambda, \nu)]$).

Given a partition of $\Lambda$ into finitely many definable pieces, it suffices to prove the proposition for each piece.
Moreover, given a definable bijection $\phi\colon \Lambda' \to \Lambda$, compatible with the projection to $S$,
Proposition~\ref{prop.Z-to-Qp} allows us to replace $P(\Lambda, \nu)$ by $P(\Lambda', \nu \circ \phi)$.
In the following, we will apply those two techniques various times to simplify $\Lambda$.

By appyling the Parametric Rectilinearization \cite[Theorem 3]{Clu.cell} to $\Lambda$, we reduce to the case
where $\Lambda$ is of the form $\Lambda_0 \times \NN^m$
for some non-empty $\Lambda_0 \sub S \times \ZZ^{n - m}$ having finite fibers $\Lambda_{0,s}$, and where the function $\nu\colon \Lambda \to \ZZ$ is of the form
\begin{equation}\label{eq.nu1}
\nu(s, \lambda_1, \dots, \lambda_n) = \frac1r (c(s) + \sum_i b_i \lambda_i)
\end{equation}
for some integers $r \ge 1$ and $b_i$ and some definable $c\colon S \to \ZZ$.
Using that $\nu$ takes integer values on all of $\Lambda$, we deduce that $b_i$ is a multiple of $r$ for each $i > n - m$.
Then, Lemma~\ref{lem.pres-prod} allows us to write $[P(\Lambda, \nu)]$ as a product
\[
[P(\Lambda_0, \nu_0)] \cdot \prod_{i=1}^m[P(S \times \NN, \nu_i)],
\]
for $\nu_i(s, \lambda_i) = \frac{b_i}{r}\lambda_i$ and $\nu_0$ the ``rest''. Thus it remains to treat the following two cases:

\medskip

Case 1. $X = P(S \times \NN, \nu)$ with $\nu(s,\lambda) = b\lambda$ for some $b \in \ZZ$.

\medskip

Since nothing depends on $s$, we omit $s$ and $S$ in the proof of this case.

The $p$-adic measure of $X$ is equal to $\sum_{\lambda \in \NN} p^{b\lambda}$, so this measure being finite implies
$b < 0$. Applying Proposition~\ref{prop.Z-to-Qp} to the diagonal embedding
\[
\Delta\colon \NN \to \ZZ^{n} , \lambda \mapsto (\lambda, \dots, \lambda)\]
for $n := -b$ yields that $[X]$ equals the class of $P(\Delta(\NN), (\lambda, \dots, \lambda) \mapsto \nu(\lambda)))$.
By definition (and using $\nu(\lambda) = -n\lambda$), this set is equal to
$$
\{(x,y)\in P(\Delta(\NN)) \times \QQ_p \mid \ac(y) = 1, v(y) = -n-1\},
$$
which is the Cartesian product of $P(\Delta(\NN))$ with a translate of $p^{-n}\ZZ_p$.
Now we conlude using Lemma~\ref{lem.pn-1}:
\[
(p^n-1)[X] = (p^n-1)[P(\Delta(\NN))] \cdot [p^{-n}\ZZ_p]
= 1 \cdot p^n \in \Rbasic.
\]

\medskip

Case 2. $X=P(\Lambda,\nu)$, where $\Lambda \sub S \times \ZZ^n$ has finite fibers over $S$ and
\begin{equation}\label{eq.nu}
\nu(s, \lambda_1, \dots, \lambda_n) = \frac1r (c(s) + \sum_i b_i \lambda_i)
\end{equation}
with $r, c, b_i$ as in (\ref{eq.nu1}).

\medskip

If all $b_i$ are $0$, then $X$ is a basic set and we are done, so suppose without loss that $b_n \ne 0$. We do an induction on $n$, i.e., we will reduce this to a case where $X = P(\Lambda', \nu')$ for some $\Lambda' \sub S \times \ZZ^{n-1}$ with finite fibers over $S$. Note that if $n = 0$, then clearly $X$ is basic.

Without loss, $b_n < 0$. By cell-decomposing $\Lambda$ with respect to the last variable, we may assume that it is of the form
\[
\Lambda = \{(\hat\lambda, \lambda_n) \in \hat\Lambda \times \ZZ
\mid \alpha_1(\hat\lambda) \le \lambda_n < \alpha_2(\hat{\lambda}),\
\lambda_n \equiv \mu \mod \ell\}
\]
for some integers $\mu \ge 0, \ell > 0$, some definable $\hat{\Lambda} \sub S \times \ZZ^{n-1}$ and some definable
$\alpha_1, \alpha_2 \colon \hat{\Lambda} \to \ZZ$ satisfying $\alpha_1(\hat \lambda)<\alpha_2(\hat \lambda)$ for all $\hat\lambda \in \hat\Lambda$.
By a further partition and a linear transformation,
we may get rid of the congruence condition, so that the set $\Lambda$ is equal to the set-theoretic difference
$\Lambda_1 \setminus \Lambda_2$, where
\[
\Lambda_i = \{(\hat\lambda, \lambda_n) \in \hat\Lambda \times \ZZ
\mid \alpha_i(\hat\lambda) \le \lambda_n\}
\]
for $i = 1, 2$. Set $X_i := P(\Lambda_i, \nu)$ (where $\nu$ is extended to the larger domain using Equation~(\ref{eq.nu})). Since $b_n < 0$, the sets $X_{i,s}$ have finite measure and we have the equation $[X] = [X_1] - [X_2]$ in $R_S$.
In this way, we reduced the problem of proving that a multiple of $[X]$ lies in $\Rbasic_S$ to proving it for both $[X_i]$.
Using a linear transformation, we reduce to the case $\alpha_i = 0$, so that $\Lambda_i = \Lambda'_i\times \NN$ (for some $\Lambda'_i \sub S\times \ZZ^{n-1}$ with finite fibers over $S$).
As in the above discussion just before Case 1, we now can write $[X_i]$ as a product
$[P(\Lambda'_i,\nu'_i)] \cdot [P(\NN,\nu''_i)]$.
The first factor is treated by the induction on $n$ in Case 2, and the second one has already been treated before, in Case 1. This finishes the proof of the proposition.
\end{proof}

\begin{proof}[Proof of Theorem~\ref{thm.main}]
Let $\Xi \in R_S$ be given such that $\mu_s(\Xi) = 0$ for every $s \in S$. (We need to prove $\Xi = 0$.)

By Proposition~\ref{prop.make.basic.fam}, we find an integer $\ell \ge 1$ such that
$\ell\cdot\Xi \in \Rbasic_S$. Then we also have $\mu_s(\ell\cdot\Xi) = 0$ for every $s \in S$,
so by Proposition~\ref{prop.basic}, we obtain
$\ell\cdot\Xi = 0$. Now $\Xi = 0$ follows from $R_S$ being torsion free (Lemma~\ref{lem.div}).
\end{proof}

\begin{rem}\label{rem:genK-L}
The only ingredients used in the entire paper are cell decomposition (which implies dimension theory and related results), the existence of definable Skolem functions and the fact that the (imaginary) $\cL$-definable subsets of $\ZZ^n$ are exactly the Presburger sets. Therefore, all our results (notably Theorem \ref{thm.main}) also hold in various generalized situations where we have those ingredients. In particular:
\begin{enumerate}
 \item Instead of $\cL$ as specified at the beginning of Section~\ref{sec.results} (as parameter free language of valued fields), one can take any language which expands $\cL$ with an analytic structure on $\QQ_p$, as in \cite{CL.analyt}. (The ingredients hold in this generality by \cite{CL.analyt}.) In particular, the map from Theorem \ref{thm.main} is still injective if we define the ring $R_S$ using subanalytic sets on $\QQ_p$, in the sense of \cite{DD.suban}. Indeed, \cite{CL.analyt}.
 \item Instead of $\QQ_p$, any finite field extension $K$ of $\QQ_p$ can be used, provided that one expands the language $\cL$ with a constant symbol for a uniformizing element $\varpi$ and enough constant symbols to obtain definable Skolem functions. (And again, one can expand $\cL$ by an analytic structure.)

 To make the proofs work in $K$, most occurences of $p$ need to be replaced either by $\varpi$ or by the cardinality of the residue field. The least straight forward changes might be those to the proof of Lemma~\ref{lem.div}: There, we obtain $[\varpi\cO_K]$ as a multiplicative inverse of $q$ (which then yields that $p$ is invertible since $p \mid q$),
 and to get that $\ell \ne p$ is invertible, we use $\ell \mid (p^n - 1) \mid (q^n - 1)$, and invertibility of $q^n-1$ follows from Lemma~\ref{lem.pn-1}.
 \item
 Even more generally, one can use any langugage $\cL$ such that the $\cL$-theory of $K$ is
 hensel minimal (more precisely, $\omega$-\heqc-minimal as defined in \cite[Section 6.1]{iCR.hmin}), with pure Presburger structure on the value group and with definable Skolem functions on $K$. (For cell decomposition, see  \cite[Theorem~5.4.2 and Addendum~1]{iCR.hmin})
\end{enumerate}
\end{rem}

\bibliographystyle{amsplain}
\bibliography{references}

\providecommand{\bysame}{\leavevmode\hbox to3em{\hrulefill}\thinspace}
\providecommand{\MR}{\relax\ifhmode\unskip\space\fi MR }
\providecommand{\MRhref}[2]{%
  \href{http://www.ams.org/mathscinet-getitem?mr=#1}{#2}
}
\providecommand{\href}[2]{#2}
\begin{thebibliography}{10}

\bibitem{CL.analyt}
R.~Cluckers and L.~Lipshitz, \emph{{Fields with analytic structure.}}, J. Eur.
  Math. Soc. (JEMS) \textbf{13} (2011), no.~4, 1147--1223 (English).

\bibitem{Clu.DefQp}
Raf Cluckers, \emph{Classification of semi-algebraic {$p$}-adic sets up to
  semi-algebraic bijection}, J. Reine Angew. Math. \textbf{540} (2001),
  105--114. \MR{MR1868600 (2002i:14052)}

\bibitem{Clu.cell}
\bysame, \emph{Presburger sets and {$p$}-minimal fields}, J. Symbolic Logic
  \textbf{68} (2003), no.~1, 153--162. \MR{MR1959315 (2003m:03062)}

\bibitem{CCL.lipschitz}
Raf Cluckers, Georges Comte, and Fran{\c{c}}ois Loeser, \emph{Lipschitz
  continuity properties for {$p$}-adic semi-algebraic and subanalytic
  functions}, Geom. Funct. Anal. \textbf{20} (2010), no.~1, 68--87.
  \MR{2647135}

\bibitem{iC.grZ}
Raf Cluckers and Immanuel Halupczok, \emph{Definable sets up to definable
  bijections in presburger groups}, Transactions of the London Mathematical
  Society \textbf{5} (2018), no.~1, 47--70.

\bibitem{iCR.hmin}
Raf Cluckers, Immanuel Halupczok, and Silvain Rideau-Kikuchi, \emph{Hensel
  minimality},  (2019).

\bibitem{CL.mot}
Raf Cluckers and Fran{\c{c}}ois Loeser, \emph{Constructible motivic functions
  and motivic integration}, Invent. Math. \textbf{173} (2008), no.~1, 23--121.

\bibitem{DD.suban}
J.~Denef and Lou van~den Dries, \emph{{$p$}-adic and real subanalytic sets},
  Ann. of Math. (2) \textbf{128} (1988), no.~1, 79--138. \MR{MR951508
  (89k:03034)}

\bibitem{Den.rat}
Jan Denef, \emph{The rationality of the {P}oincar\'e series associated to the
  {$p$}-adic points on a variety}, Invent. Math. \textbf{77} (1984), no.~1,
  1--23. \MR{MR751129 (86c:11043)}

\bibitem{Den.cell}
\bysame, \emph{{$p$}-adic semi-algebraic sets and cell decomposition}, J. Reine
  Angew. Math. \textbf{369} (1986), 154--166. \MR{MR850632 (88d:11030)}

\bibitem{HK.motInt}
Ehud Hrushovski and David Kazhdan, \emph{Integration in valued fields},
  Algebraic geometry and number theory, Progr. Math., vol. 253, Birkh\"auser
  Boston, Boston, MA, 2006, pp.~261--405. \MR{MR2263194 (2007k:03094)}

\bibitem{KZ.per}
Maxim Kontsevich and Don Zagier, \emph{Periods}, Mathematics unlimited---2001
  and beyond, Springer, Berlin, 2001, pp.~771--808. \MR{1852188 (2002i:11002)}

\bibitem{Mac.QE}
Angus Macintyre, \emph{On definable subsets of {$p$}-adic fields}, J. Symbolic
  Logic \textbf{41} (1976), no.~3, 605--610. \MR{MR0485335 (58 \#5182)}

\bibitem{Pre.oagZ}
M.~Presburger, \emph{{\"Uber die Vollst\"andigkeit eines gewissen Systems der
  Arithmetik ganzer Zahlen, in welchem die Addition als einzige Operation
  hervortritt.}}, {}, 1930 (German).

\bibitem{Dri.skolQp}
Lou van~den Dries, \emph{Algebraic theories with definable {S}kolem functions},
  J. Symbolic Logic \textbf{49} (1984), no.~2, 625--629. \MR{745390}

\bibitem{Dri.dimDef}
\bysame, \emph{Dimension of definable sets, algebraic boundedness and
  {H}enselian fields}, Ann. Pure Appl. Logic \textbf{45} (1989), no.~2,
  189--209, Stability in model theory, II (Trento, 1987). \MR{MR1044124
  (91k:03082)}

\end{thebibliography}

\end{document}